\documentclass[a4paper,10pt]{article}
\usepackage{latexsym,amssymb,amsmath,amsthm,amsxtra,xypic}
\usepackage{amsfonts}
\usepackage{amscd}
\usepackage{mathrsfs}
\usepackage[dvips]{graphicx}
\usepackage[all]{xy}
\usepackage{ulem}

\newcommand{\allowpagebreak}
\allowdisplaybreaks

\newtheorem{thm}{Theorem}[section]
\newtheorem{lem}[thm]{Lemma}
\newtheorem{cor}[thm]{Corollary}
\newtheorem{pro}[thm]{Proposition}

\theoremstyle{definition}
\newtheorem{defi}[thm]{Definition}

\newcommand {\emptycomment}[1]{}

\newcommand{\pf}{\noindent{\bf Proof.}\ }

\def\DD{D}

\newcommand{\al}{\alpha}
\newcommand{\be}{\beta}

\newcommand{\bebe}{\beta^2}

\newcommand{\lam}{\lambda}
\newcommand{\si}{\sigma}

\newcommand{\h}{V}

\newcommand{\frkh}{V}

\newcommand{\dlam}{_\lambda}

\newcommand{\hal}{\hat{\alpha}}

\newcommand{\E}{\mathrm{E}}

\newcommand{\Ker}{\mathrm{Ker}}

\newcommand{\End}{\mathrm{End}}
\newcommand{\Ext}{\mathrm{Ext}}

\newcommand{\id}{\mathrm{id}}

\newcommand{\bea}{\begin{eqnarray}}
\newcommand{\eea}{\end{eqnarray}}
\newcommand{\Bea}{\begin{eqnarray*}}
\newcommand{\Eea}{\end{eqnarray*}}
\newcommand{\beq}{\begin{equation}}
\newcommand{\eeq}{\end{equation}}
\newcommand{\Beq}{\begin{equation*}}
\newcommand{\Eeq}{\end{equation*}}

\begin{document}

\allowdisplaybreaks

\title{Representations and cohomologies of Hom-Lie-Yamaguti algebras with applications}

\author{Tao Zhang$^{a}$\thanks{Corresponding author. {E-mail address:} zhangtao@htu.cn}\\
{\footnotesize $^a$ College of Mathematics and Information Science, Henan Normal University, Xinxiang 453007, PR China}}

\date{}
\maketitle



\begin{abstract}
  The representation and cohomology theory of Hom-Lie-Yamaguti algebras is introduced.
  As an application, we study deformation and extension of Hom-Lie-Yamaguti algebras.
  It proved that a 1-parameter infinitesimal deformation of a Hom-Lie-Yamaguti algebra $T$
  corresponds to a Hom-Lie-Yamaguti algebra of deformation type and a (2,3)-cocycle of $T$ with coefficients in the adjoint representation.
  We also prove that abelian extensions of Hom-Lie-Yamaguti algebras are classified by the (2,3)-cohomology group.
\end{abstract}

\section{Introduction}

In recent years, Hom-type algebras were studied by many researchers.
The first examples coming from $q$-deformations of Witt and Virasora algebras are Hom-Lie algebras, see \cite{Hom1}.
Other types include Hom-associative algebras, Hom-Nambu-Lie algebras, Hom-Hopf algebras, etc.
See \cite{Makhlouf01,Makhlouf02,Makhlouf03,Caenepeel,Sheng,Yau1,Yau2} and the references therein.

In \cite{GNI}, the authors introduced the concept of Hom-Lie-Yamaguti algebras.
It is a Hom-type generalization of a Lie-Yamaguti algebra in \cite{KW,BEM01}, a general Lie triple system in \cite{Ya1,Ya2}
and a Lie triple algebra in \cite{Kik01}.
In \cite{Ma}, the authors studied the formal deformations of Hom-Lie-Yamaguti algebras,
where only low dimensional deformation cohomology were defined without the help of any representation.
So we wonder if there is a systematic study of Hom-Lie-Yamaguti algebras to give its general representation and cohomology theory?
The present paper is denoted to solve this problem.

The method of this paper is based on our recent work \cite{Zhang0,Zhang1}.
In \cite{Zhang0}, we give a new characterization of the representation and cohomology theory of Lie triple systems.
In \cite{Zhang1}, we give a detailed study on the $(2,3)$-cohomology group associated to a representation of a Lie-Yamaguti algebra.
As an application, we study deformation and extension theory of Lie-Yamaguti algebras.
In this paper, we will first define the representation and  cohomology theory of Hom-Lie-Yamaguti algebras.
Then we will study deformation and extension theory of Hom-Lie-Yamaguti algebras as we did in \cite{Zhang1}.
We will see that they are classified by the $(2,3)$-cohomology groups.
The difficulty in this case is that, we have a morphism $\al:T\to T$ and more conditions with respect to $\al$ to be
compatible with the structure of Hom-Lie-Yamaguti algebras.
Fortunately, we overcome this difficulty by using an equivalent characterization of representation and a careful analysis on the coboundary operator.
All our results in this paper will generalize Yamaguti's representation and  cohomology theory for Lie-Yamaguti algebras in \cite{Ya1,Ya2}.

The paper is organized as follows. In section 2, we introduce the concept of  representations of a Hom-Lie-Yamaguti algebra.
Then we define the coboundary operator on the cochain complex of  a Hom-Lie-Yamaguti algebra with coefficient in a representation $V$ to produce the cohomology group.
We pay special attentions to the (2,3)-cohomology group since it will used in the following context.
In Section 3, we study the infinitesimal deformation theory of Hom-Lie-Yamaguti algebras. We prove that
there is a Hom-Lie-Yamaguti algebra of deformation type
and a (2,3)-cocycle of $T$ with coefficients in the adjoint representation associated to a deformation.
We also introduce the  notion of Nijenhuis operators to describe trivial deformations.
In Section 4, we study  abelian extensions of Hom-Lie-Yamaguti algebras.
We prove that there is a one-to-one correspondence between equivalence classes of
abelian extensions of the Hom-Lie-Yamaguti algebra $T$ by $V$ and elements of the (2,3)-cohomology group.

Throughout this paper, we work on an algebraically closed field $\mathbb{K}$ of characteristic different from 2 and 3.

\section{Representations and Cohomologies}\label{sec:LYAcohomology}

In this section, we first recall some basic definitions regarding Hom-Lie-Yamaguti algebras.
Then we define the representation and cohomology theory of Hom-Lie-Yamaguti algebras.

A Hom-vector space is a pair $(V,\al)$ in which $V$ is a vector space and $\al: V \to V$ is a linear map. A morphism
$(V,\al)\to (W,\be) $ of Hom--vector space is a linear map $f : V \to W$ such that $\be\circ f = f \circ \al$. We will often
abbreviate a Hom-vector space $(V,\al)$ to $V$.

\begin{defi}\label{def:LYA}
A Hom-Lie-Yamaguti algebra (or HLYA for short) consists of a vector space $T$ together with a linear map $\al:T\to T$,
a bilinear map $[\cdot,\cdot]:  T\times T\to  T$ and a trilinear map
$[\cdot,\cdot,\cdot]:  T\times T\times T\to  T$ such that, for all $ x_i, y_i \in  T$, the following conditions are satisfied:
\begin{itemize}
\item[$\bullet$] {\rm(HLY01)}\quad $\al([x_1, x_2])=[\al (x_1),\al (x_2)];$
\item[$\bullet$] {\rm(HLY02)}\quad $\al([x_1, x_2,x_3])=[\al (x_1), \al (x_2),\al (x_3)];$
\item[$\bullet$] {\rm(HLY1)}\quad $[x_1,x_2]+[x_2,x_1]=0$;
\item[$\bullet$] {\rm(HLY2)}\quad $[x_1, x_2,x_3]+[x_2, x_1,x_3]=0$;
\item[$\bullet$] {\rm(HLY3)}\quad $[[x_1, x_2],\al(x_3)]+c.p.+[x_1, x_2,x_3]+c.p.=0$;
\item[$\bullet$] {\rm(HLY4)}\quad $[[x_1, x_2],\al(x_3),\al(y_1)]+[[x_2, x_3],\al(x_1),\al(y_1)]+[[x_3, x_1],\al(x_2),\al(y_1)]=0$;
\item[$\bullet$] {\rm(HLY5)}\quad $[\al(x_1), \al(x_2), [y_1, y_2]] = [[ x_1, x_2, y_1], \al^2(y_2)] + [\al^2(y_1), [ x_1, x_2, y_2]]$;
\item[$\bullet$] {\rm(HLY6)}\quad
  $[\al^2(x_1), \al^2(x_2), [y_1, y_2, y_3]] = [[ x_1, x_2, y_1], \al^2(y_2), \al^2(y_3)] + [\al^2(y_1), [ x_1, x_2, y_2], \al^2(y_3)]
  + [\al^2(y_1), \al^2(y_2), [ x_1, x_2, y_3]]$.
\end{itemize}
where c.p. means cyclic permutations with respect to $x_1, x_2, x_3$.
We denote a HLYA by $(T,[\cdot,\cdot],[\cdot,\cdot,\cdot],\al)$ or simply by $T$.
\end{defi}

A linear map $\al$ satisfying (HLY01) and (HLY02) is called an algebraic homomorphism.
When $\al=\id$, conditions (HLY01) and (HLY02) are trivial and the other conditions (HLY1)--(HLY6)
reduced to  conditions (LY1)--(LY6) for a Lie-Yamaguti algebra (LYA for short) in  \cite{Zhang1}.
Note that conditions (HLY1) and (HLY2) are equivalent to
$[x_1,x_1]=0 \  \mbox{and}\  [x_1,x_1,x_3]=0$
respectively.

A homomorphism between two HLYAs $T$ and $T'$ is a linear map $\varphi:T\to T'$ satisfying $\varphi\circ \al=\al'\circ \varphi$ and
\begin{eqnarray}
\varphi([x_1, x_2])=[\varphi (x_1), \varphi (x_2)]',\quad \varphi([x_1, x_2,x_3])=[\varphi (x_1), \varphi (x_2),\varphi (x_3)]'.
\end{eqnarray}

\begin{defi}
A HLYA of {\bf deformation type} consists of a vector space $T$ together with a linear map $\al:T\to T$,
a bilinear map $\nu:  T\times T\to  T$ and a trilinear map
$\omega:  T\times T\times T\to  T$ satisfying all conditions
in Definition \ref{def:LYA} except that (HLY3) is replaced by the following conditions:
\begin{itemize}
\item[$\bullet$] {\rm(HLY3')}\quad $\nu(\nu(x_1, x_2),\al(x_3))+c.p.=0$.
\end{itemize}
\end{defi}

Now we give the definition of a representation of a HLYA.

\begin{defi}Let $(T,\al)$ be a HLYA and $(V,\be)$ be a Hom-vector space.
A representation of $(T,\al)$ on $(V,\be)$ consists of a linear map $\rho$: $T \to \End(V)$ and bilinear maps
$D, \theta$: $T\times T\to \End(V)$ such that the following conditions are satisfied:
\begin{itemize}
\item[$\bullet$]{\rm(HR01)}\quad $\rho(\al(x_1))\circ \be=\be\circ\rho(x_1)$;
\item[$\bullet$]{\rm(HR02)}\quad $\DD(\al(x_1),\al(x_2))\circ \be=\be\circ\DD(x_1,x_2)$;
\item[$\bullet$]{\rm(HR03)}\quad $\theta(\al(x_1),\al(x_2))\circ \be=\be\circ\theta(x_1,x_2)$;
\item[$\bullet$]{\rm(HR31)}\quad $\DD(x_1,x_2)-\theta(x_2,x_1)+\theta(x_1,x_2)+\rho([x_1,x_2])\circ\be-\rho( \al(x_1))\rho (x_2)+\rho (\al(x_2))\rho (x_1)=0$;
\item[$\bullet$]{\rm(HR41)}\quad $\DD([x_1,x_2],\al(x_3))+\DD([x_2,x_3],\al(x_1))+\DD([x_3,x_1],\al(x_2))=0$;
\item[$\bullet$]{\rm(HR42)}\quad $\theta([x_1,x_2],\al(y_1))\circ \be=\theta(\al(x_1),\al(y_1))\rho(x_2)-\theta(\al(x_2),\al(y_1))\rho(x_1)$;
\item[$\bullet$]{\rm(HR51)}\quad $\DD(\al(x_1),\al(x_2))\rho(y_2)=\rho(\al^2(y_2))\DD(x_1,x_2)+\rho([x_1, x_2, y_2])\circ\bebe.$;
\item[$\bullet$]{\rm(HR52)}\quad $\theta(\al(x_1),[y_1, y_2])\circ \be=\rho(\al^2(y_1))\theta(x_1, y_2)-\rho(\al^2(y_2))\theta(x_1,y_1)$;
\item[$\bullet$]{\rm(HR61)}\quad  $\DD(\al^2(x_1),\al^2(x_2))\theta(y_1,y_2)$
     \begin{eqnarray*}
     &=&\theta(\al^2(y_1),\al^2(y_2))\DD(x_1,x_2)
           +\theta([x_1,x_2,y_1],\al^2(y_2))\circ\be^2+\theta(\al^2(y_1),[x_1,x_2,y_2])\circ\be^2;
       \end{eqnarray*}
\item[$\bullet$]{\rm(HR62)}\quad  $\theta(\al^2(x_1),[y_1, y_2, y_3])\circ\be^2$
   \begin{eqnarray*}
    &=& \theta (\al^2(y_2), \al^2(y_3))\theta(x_1,y_1)
    -\theta (\al^2(y_1), \al^2(y_3))\theta(x_1,y_2) + \DD (\al^2(y_1), \al^2(y_2))\theta(x_1,y_3).
   \end{eqnarray*}
\end{itemize}
In this case, we also call $V$ to be a $T$-module.
\end{defi}

For example, given a HLYA $T$, there is a natural adjoint representation on itself.
The corresponding representation maps $\rho, D$ and $\theta$ are given by
\begin{eqnarray*}
 \rho(x_1)(x_2):=[x_1,x_2],\quad D(x_1,x_2)x_3:=[x_1,x_2,x_3],\quad \theta(x_1,x_2)x_3:=[x_3,x_1,x_2].
\end{eqnarray*}

The next Proposition \ref{prop:33} gives an equivalent characterization of a representation.
The proof is omitted since it is same as the proof of Lemma \ref{lem:33} in the last section.
\begin{pro}\label{prop:33}
Let $(T,\al)$ be a  HLYA and $(V,\be)$ be a Hom-vector space.
Assume we have a map $\rho$ from $T $ to $\End(V)$ and maps $D, \theta$ from $T\times T$ to $\End(V)$.
Then $(\rho,\DD,\theta)$  is a representation of $T$ on $V$ if and only if $T\oplus V$ is a HLYA under the following maps:
\begin{eqnarray}
\label{semi:newbracket00}{(\al+\be)(x_1 + u_1)}&\triangleq&\al(x_1)+\be(u_1),\\
\label{semi:newbracket01}{[x_1 + u_1, x_2 + u_2]}&\triangleq&[x_1, x_2] +\rho(x_1)(u_2) -\rho(x_2)(u_1),\\
\notag{[x_1 + u_1, x_2 + u_2, x_3 + u_3]}&\triangleq&[x_1, x_2, x_3] + \DD(x_1, x_2)(u_3)-\theta(x_1, x_3)(u_2)\\
\label{semi:newbracket02}&& + \theta( x_2, x_3)(u_1),
\end{eqnarray}
In this case, the HLYA $T\oplus V$ is called semidirect product of $T$ and $V$,
which is denoted by $T\ltimes V$.
\end{pro}

\medskip

Motivated by Yamaguti's cohomology for Lie-Yamaguti algebras, we are going to define cohomology for HLYAs as follows.

Let $V$ be a representation of HLYA $T$. Let us define the cohomology groups of $T$ with coefficients in $V$.
Let $f:T\times\cdots \times T\to V$ be $n$-linear maps of $T$ into $V$ such that the following conditions are satisfied:
\begin{eqnarray}
\label{eq:cochain01} f(\al(x_1)\cdots,\al(x_n))&=&\be(f(x_1,\cdots,x_n)),\\
 f(x_1,\cdots,x_{2i-1},x_{2i}\cdots,x_n)&=&0,\ \  \mbox{if}\ \  x_{2i-1}=x_{2i}.
\end{eqnarray}
The vector space spanned by such linear maps is called an $n$-cochain of $T$, which is denoted by $C^n(T,V)$ for $n\geq 1$.

\begin{defi}\label{def:cohomology}
For any $(f,g)\in C^{2n}(T, V)\times C^{2n+1}(T, V)$
the coboundary operator $\delta: (f, g)\mapsto (\delta_{\textrm{I}}f, \delta_{\textrm{II}}g)$
is a mapping from $C^{2n}(T, V)\times C^{2n+1}(T, V)$ into $C^{2n+2}(T, V)\times C^{2n+3}(T, V)$ defined as follows:
\begin{eqnarray*}
 &&(\delta_{\textrm{I}}f)(x_{1}, x_{2}, \cdots, x_{2n+2})\\
&=&\rho(\al^{2n}(x_{2n+1}))g(x_{1}, \cdots, x_{2n}, x_{2n+2}))-\rho(\al^{2n}(x_{2n+2}))g(x_{1}, \cdots, x_{2n+1})\\
&&-g(\al(x_{1}), \cdots, \al(x_{2n}), [x_{2n+1},x_{2n+2}])\\
&&+\sum\limits_{k=1}^{n}(-1)^{n+k+1}D(\al^{2n-1}(x_{2k-1}), \al^{2n-1}(x_{2k}))f(x_{1}, \cdots, \hat{x}_{2k-1}, \hat{x}_{2k}, \cdots, x_{2n+2})\\
&&+\sum\limits_{k=1}^{n}\sum\limits_{j=2k+1}^{2n+2}(-1)^{n+k}f(\al^2(x_{1}), \cdots, \hat{x}_{2k-1}, \hat{x}_{2k}, \cdots, [x_{2k-1}, x_{2k}, x_{j}], \cdots, \al^2(x_{2n+2})),
\end{eqnarray*}
\begin{eqnarray*}
 &&(\delta_{\textrm{II}}g)(x_{1}, x_{2}, \cdots, x_{2n+3})\\
&=& \theta(\al^{2n}(x_{2n+2}), \al^{2n}(x_{2n+3}))g(x_{1},\cdots, x_{2n+1})\\
&&-\theta(\al^{2n}(x_{2n+1}), \al^{2n}(x_{2n+3}))g(x_{1}, \cdots, x_{2n}, x_{2n+2})\\
&&+\sum\limits_{k=1}^{n+1}(-1)^{n+k+1}D(\al^{2n}(x_{2k-1}), \al^{2n}(x_{2k}))g(x_{1}, \cdots, \hat{x}_{2k-1}, \hat{x}_{2k}, \cdots, x_{2n+3})\\
&&+\sum\limits_{k=1}^{n+1}\sum\limits_{j=2k+1}^{2n+3}(-1)^{n+k}g(\al^2(x_{1}), \cdots, \hat{x}_{2k-1}, \hat{x}_{2k}, \cdots, [x_{2k-1}, x_{2k}, x_{j}], \cdots, \al^2(x_{2n+3})).
\end{eqnarray*}
\end{defi}

When $\al=\id$, one recovers Yamaguti's cohomology for LYA in \cite{Ya2}.

\begin{lem}
With the above notations, for any $(f,g)\in C^{2n}(T, V)\times C^{2n+1}(T, V)$, we have
\begin{eqnarray}
\label{eq001}\delta_{\textrm{I}} f(\al(x_1)\cdots,\al(x_{2n+2}))&=&\be(\delta_{\textrm{I}} f(x_{1}, x_{2}, \cdots, x_{2n+2})),\\
\label{eq002}\delta_{\textrm{II}} g(\al(x_1)\cdots,\al(x_{2n+3}))&=&\be(\delta_{\textrm{II}} g(x_{1}, x_{2}, \cdots, x_{2n+3})).
\end{eqnarray}
Thus we obtain a well-defined map
 $$\delta=(\delta_{\textrm{I}}, \delta_{\textrm{II}}): C^{2n}(T, V)\times C^{2n+1}(T, V)\to C^{2n+2}(T, V)\times C^{2n+3}(T, V).$$
\end{lem}

\pf We only prove equation \eqref{eq001} since equation \eqref{eq002} can be verified similarly. By Definition \ref{def:cohomology}, we have
\begin{eqnarray*}
 &&(\delta_{\textrm{I}}f)(\al(x_1),\cdots,\al(x_{2n+2}))\\
&=&\rho(\al^{2n+1}(x_{2n+1}))g(\al(x_1), \cdots, \al(x_{2n+2}))-\rho(\al^{2n+1}(x_{2n+2}))g(\al(x_1), \cdots, \al(x_{2n+1}))\\
&&-g(\al^2(x_{1}), \cdots, \al^2(x_{2n}), \al([x_{2n+1},x_{2n+2}]))\\
&&+\sum\limits_{k=1}^{n}(-1)^{n+k+1}D(\al^{2n}(x_{2k-1}), \al^{2n}(x_{2k}))f(\al(x_1), \cdots, \hat{x}_{2k-1}, \hat{x}_{2k}, \cdots, \al(x_{2n+2}))\\
&&+\sum\limits_{k=1}^{n}\sum\limits_{j=2k+1}^{2n+2}(-1)^{n+k}f(\al^3(x_1), \cdots, \hat{x}_{2k-1}, \hat{x}_{2k}, \cdots,
\al([x_{2k-1}, x_{2k}, x_{j}]), \cdots, \al^3(x_{2n+2}))\\
&=&\rho(\al^{2n+1}(x_{2n+1}))\circ\be(g(x_1, \cdots, x_{2n+2}))-\rho(\al^{2n+1}(x_{2n+2}))\circ\be( g(x_1, \cdots, x_{2n+1}))\\
&&-\be(g(\al(x_{1}), \cdots, \al(x_{2n}), [x_{2n+1},x_{2n+2}]))\\
&&+\sum\limits_{k=1}^{n}(-1)^{n+k+1}D(\al^{2n}(x_{2k-1}), \al^{2n}(x_{2k}))\circ\be( f(x_1, \cdots, \hat{x}_{2k-1}, \hat{x}_{2k}, \cdots, x_{2n+2}))\\
&&+\sum\limits_{k=1}^{n}\sum\limits_{j=2k+1}^{2n+2}(-1)^{n+k}\be\Big (f(\al^2(x_1), \cdots, \hat{x}_{2k-1}, \hat{x}_{2k}, \cdots,
[x_{2k-1}, x_{2k}, x_{j}], \cdots, \al^2(x_{2n+2}))\Big)\\
&=&\be(\rho(\al^{2n}(x_{2n+1}))g(x_{1}, \cdots, x_{2n}, x_{2n+2})))-\be(\rho(\al^{2n}(x_{2n+2}))g(x_{1}, \cdots, x_{2n+1}))\\
&&-\be( g(\al(x_{1}), \cdots, \al(x_{2n}), [x_{2n+1},x_{2n+2}]))\\
&&+\sum\limits_{k=1}^{n}(-1)^{n+k+1}\be(D(\al^{2n-1}(x_{2k-1}), \al^{2n-1}(x_{2k}))f(x_{1}, \cdots, \hat{x}_{2k-1}, \hat{x}_{2k}, \cdots, x_{2n+2}))\\
&&+\sum\limits_{k=1}^{n}\sum\limits_{j=2k+1}^{2n+2}(-1)^{n+k}\be\Big(f(\al^2(x_{1}), \cdots, \hat{x}_{2k-1}, \hat{x}_{2k}, \cdots, [x_{2k-1}, x_{2k}, x_{j}], \cdots, \al^2(x_{2n+2}))\Big)\\
&=&\be((\delta_{\textrm{I}}f)(x_{1}, x_{2}, \cdots, x_{2n+2}))
\end{eqnarray*}
where in the second equality we use condition \eqref{eq:cochain01} and in the third equality we use conditions (HR01)--(HR03).
\qed

\begin{pro}\label{pro:coboundary}
The coboundary operator defined above satisfies
$\delta\circ \delta=0$, that is $\delta_{\textrm{I}} \circ\delta_{\textrm{I}} =0$ and $\delta_{\textrm{II}} \circ\delta_{\textrm{II}} =0.$
\end{pro}

The above  Proposition \ref{pro:coboundary} can be proved by tedious computations.
For our purpose, we will check a special case in Proposition \ref{prop:1cob}.
\medskip

Let $Z^{2n}(T, V)\times Z^{2n+1}(T, V)$ be the subspace of $C^{2n}(T, V)\times C^{2n+1}(T, V)$ spanned by $(f, g)$ such that $\delta(f, g)=0$
which is called the space of cocycles
and $B^{2n}(T, V)\times B^{2n+1}(T, V)=\delta(C^{2n-2}(T, V)\times C^{2n-1}(T, V))$ which is called the space of coboundaries.

\begin{defi}\label{def:1cohmologygroup} For the case $n\geq 2$,
the $(2n, 2n+1)$-cohomology group  of a HLYA $T$ with coefficients in $V$ is defined to be the quotient space:
$$H^{2n}(T, V)\times H^{2n+1}(T, V)\triangleq(Z^{2n}(T, V)\times Z^{2n+1}(T, V))/(B^{2n}(T, V)\times B^{2n+1}(T, V)).$$
\end{defi}

In conclusion, we obtain a cochain complex whose cohomology group is called cohomology group of a HLYA $T$ with coefficients in $V$.

\medskip

Note that in the above Definition \ref{def:1cohmologygroup} we assume $n\geq 2$.
For the case $n=1$, we define the (2,3)-cohomology group of a HLYA $T$ with coefficients in $V$ as follows.

\medskip

Let $C^2(T,V)$ be the space of maps $\nu:  T\times T\to \h$ such that $\nu(x_1,x_2)=-\nu(x_2,x_1)$ and
\begin{itemize}
\item[$\bullet$]{\rm(CC01)}\quad $\nu(\al(x_1),\al(x_2))=\be\circ \nu(x_1,x_2)$.
\end{itemize}
Let $C^3(T,V)$ be the space of maps $\omega:  T\times T\times T\to \h$ such that $\omega(x_1, x_2, x_3)=-\omega(x_2, x_1, x_3)$ and
\begin{itemize}
\item[$\bullet$]{\rm(CC02)}\quad $\omega(\al(x_1),\al(x_2),\al(x_3))=\be\circ \omega(x_1, x_2, x_3)$.
\end{itemize}

\begin{defi}\label{def:1coc}
Let $(T,\al)$ be a HLYA and $(V, \be)$ a $T$-module.
Then $(\nu,\omega)\in C^2(T,V)\times C^3(T,V)$ is called a (2,3)-cocycle if for all $x_1, x_2,y_1, y_2, y_3\in  T$, we have
\begin{itemize}
\item[$\bullet$]{\rm(CC1)}\quad $\omega(x_1,x_2, x_3)+c.p.-\rho(\al(x_1))\nu(x_2,x_3)-c.p.+\nu([x_1,x_2],\al(x_3))+c.p.=0$;
\item[$\bullet$]{\rm(CC2)}\quad $\theta(\al(x_1), \al(y_1))\nu(x_2,x_3)+c.p.+\omega([x_1,x_2], \al(x_3), \al(y_1))+c.p.=0$;
\item[$\bullet$]{\rm(CC3)}\quad \begin{eqnarray*}
&& \omega(\al(x_1), \al(x_2),[y_1, y_2])+\DD(\al(x_1), \al(x_2))\omega(y_1, y_2)\\
&=&\nu([x_1, x_2, y_1], \al^2(y_2)) + \nu(\al^2(y_1),[x_1, x_2,y_2]))\\
&&+\rho(\al^2(y_1))\omega(x_1,x_2,y_2)-\rho(\al^2(y_2))\omega(x_1,x_2,y_1);
\end{eqnarray*}
\item[$\bullet$]{\rm(CC4)}\quad \begin{eqnarray*}
&& \omega(\al^2(x_1), \al^2(x_2),[y_1, y_2, y_3])+\DD(\al^2(x_1), \al^2(x_2))\omega(y_1, y_2, y_3)\\
&=&\omega([x_1, x_2, y_1], \al^2(y_2), \al^2(y_3)) + \omega(\al^2(y_1),[x_1, x_2,y_2], \al^2(y_3)) \\
&&+ \omega(\al^2(y_1), \al^2(y_2), [x_1, x_2,y_3])+\theta(\al^2(y_2), \al^2(y_3))\omega(x_1,x_2,y_1)\\
&&- \theta(\al^2(y_1), \al^2(y_3))\omega(x_1,x_2,y_2) + \DD(\al^2(y_1), \al^2(y_2))\omega(x_1,x_2,y_3).
\end{eqnarray*}
\end{itemize}
The space of (2,3)-cocycles is denoted by $Z^{2}(T, V)\times Z^{3}(T, V)$.
\end{defi}
We remark that the conditions (CC3) and (CC4) are equivalent to $\delta_{\textrm{I}}(\nu)=0$ and  $\delta_{\textrm{II}}(\omega)=0$ respectively.
Why we add conditions (CC1) and (CC2) can be seen from the following context.

Let $f$ be a linear mapping of $T$ into a representation space $V$. Then $f$ is called a {\bf derivation} of $T$ into $V$ if
\begin{eqnarray}
f([x_1,x_2])&=&\rho(x_1)f(x_2)-\rho(x_2)f(x_1),\\
f([x_1, x_2, x_3])&=&\theta(x_2, x_3)f(x_1)-\theta(x_1, x_3)f(x_2)+D(x_1, x_2)f(x_3).
\end{eqnarray}

\begin{defi}\label{def:1cob}
Let $(T,\al)$ be a HLYA and $(V,\be)$ a $T$-module.
Then $(\nu,\omega)\in C^2(T,V)\times C^3(T,V)$ is called a (2,3)-coboundary if there exists a map $f: T\to V$ such that
\begin{itemize}
\item[$\bullet$]{\rm(BB01)}\quad $f\circ \al=\be\circ f$;
\item[$\bullet$]{\rm(BB1)}\quad $\nu(x_1,x_2)=\rho(x_1)f(x_2)-\rho(x_2)f(x_1)-f([x_1,x_2])$;
\item[$\bullet$]{\rm(BB2)}\quad $\omega(x_1,x_2,x_3)=\theta(x_2, x_3)f(x_1)-\theta(x_1, x_3)f(x_2)+D(x_1, x_2)f(x_3)-f([x_1, x_2, x_3])$.
\end{itemize}
The space of (2,3)-coboundaries is denoted by $B^{2}(T, V)\times B^{3}(T, V)$.
\end{defi}

\begin{pro}\label{prop:1cob}
The space of (2,3)-coboundaries is contained in space of (2,3)-cocycles.
\end{pro}

\pf
We will verify that if $(\nu,\omega)$ satisfies (BB01), (BB1) and (BB2), then it must satisfies conditions (CC01),(CC02) and (CC1)--(CC4).

By definition, for (CC01), we have
\begin{eqnarray*}
&&\nu(\al(x_1),\al(x_2))-\be\circ \nu(x_1,x_2)\\
&=& \rho(\al(x_1))f(\al(x_2))-\rho(\al(x_2))f(\al(x_1))-f([\al(x_1),\al(x_2)])\\
&&-\be\circ\{\rho(x_1)f(x_2)-\rho(x_2)f(x_1)- f([x_1,x_2])\}\\
&=& \underline{\rho(\al(x_1))\circ \be}\circ f(x_2)-\underline{\underline{\rho(\al(x_2))\circ \be}}\circ f(x_1)-{f\circ \al}([x_1,x_2])\\
&&-\underline{\be\circ\rho(x_1)}\circ f(x_2)+\underline{\underline{\be\circ\rho(x_2)}}\circ f(x_1)+{\be\circ f}([x_1,x_2])\\
&=&0.
\end{eqnarray*}
where in the last equality we have used (HR01) and (BB01).

For (CC02), we have
\begin{eqnarray*}
&&\omega(\al(x_1),\al(x_2),\al(x_3))-\be\circ \omega(x_1, x_2, x_3)\\
&=&\theta(\al(x_2), \al(x_3))f(\al(x_1))-\theta(\al(x_1), \al(x_3))f(\al(x_2))\\
&&+D(\al(x_1), \al(x_2))f(\al(x_3))-f([\al(x_1), \al(x_2), \al(x_3)])\\
&&\be\circ\{\theta(x_2, x_3)f(x_1)-\theta(x_1, x_3)f(x_2)+ D(x_1, x_2)f(x_3)- f([x_1, x_2, x_3])\}\\
&=&\underline{\theta(\al(x_2), \al(x_3))\circ \be}\circ f(x_1)-\underline{\underline{\theta(\al(x_1), \al(x_3))\circ \be}}\circ f(x_2)\\
&&+\underline{\underline{\underline{D(\al(x_1), \al(x_2))\circ \be}}}\circ f(x_3)-{f\circ \al}([x_1, x_2, x_3])\\
&&\underline{\be\circ\theta(x_2, x_3)}\circ f(x_1)-\underline{\underline{\be\circ\theta(x_1, x_3)}}\circ f(x_2)
+\underline{\underline{\underline{\be\circ D(x_1, x_2)}}}\circ f(x_3)\\
&&-{\be\circ f}([x_1, x_2, x_3])\\
&=&0.
\end{eqnarray*}
where in the last equality we have used (HR02), (HR03) and (BB01).

For (CC1), we have
\begin{eqnarray*}
&&\omega(x_1,x_2, x_3)+c.p.-\rho(\al(x_1))\nu(x_2,x_3)-c.p.+\nu([x_1,x_2],\al(x_3))+c.p.\\
&=&\Big(\theta(x_2, x_3)f(x_1)-\theta(x_1, x_3)f(x_2)+D(x_1, x_2)f(x_3)-f([x_1, x_2, x_3])\Big)+c.p.\\
&&-\rho(\al(x_1))\Big(\rho(x_2)f(x_3)-\rho(x_3)f(x_2)-f([x_2,x_3])\Big)+c.p.\\
&&+\Big(\rho([x_1,x_2])f(\al(x_3))-\rho(\al(x_3))f([x_1,x_2])-f([[x_1,x_2],\al(x_3)])\Big)+c.p.\\
&=&\Big(\DD(x_1,x_2)-\theta(x_2,x_1)+\theta(x_1,x_2)+\rho[x_1,x_2]\circ\be\\
&&-\rho (\al(x_1))\rho (x_2)+\rho (\al(x_2))\rho (x_1)\Big)f(x_3)+c.p.\\
&&-f\Big([x_1, x_2, x_3]+c.p.+[[x_1,x_2],\al(x_3)]+c.p.\Big)\\
&=&0.
\end{eqnarray*}
The last equality is by (HR31) and (HLY3).

By direct computations, for (CC2), we get
\begin{eqnarray*}
&&\theta(\al(x_1), \al(y_1))\nu(x_2,x_3)+c.p.+\omega([x_1,x_2],\al(x_3), \al(y_1))+c.p.\\
&=&\theta(\al(x_1), \al(y_1))\Big(\rho(x_2)f(x_3)-\rho(x_3)f(x_2)-f([x_2,x_3])\Big)+c.p.\\
&&+\Big(\theta(\al(x_3), \al(y_1))f([x_1,x_2])-\theta([x_1,x_2], \al(y_1))f(\al(x_3))\\
&&+D([x_1,x_2], \al(x_3))f(\al(y_1))-f([[x_1,x_2], \al(x_3), \al(y_1)])\Big)+c.p.\\
&=&\Big(\DD([x_1,x_2],\al(x_3))+\DD([x_2,x_3],\al(x_1))+\DD([x_3,x_1],\al(x_2))\Big)f(\al(y_1))\\
&&-\Big(\theta([x_1,x_2],y_1)\circ \be-\theta(\al(x_1),\al(y_1))\rho(x_2)+\theta(\al(x_2),\al(y_1))\rho(x_1)\Big)f(x_3)-c.p.\\
&&-f([[x_1,x_2], \al(x_3), \al(y_1)]+c.p.)\\
&=&0.
\end{eqnarray*}
The last equality is by (HR41), (HR42) and (HLY4).

The other cases can checked as follows:  (CC3) is valid by conditions (HR51), (HR52) and (HLY5); (CC4) is valid by conditions (HR61), (HR62) and (HLY6).
Therefore the space of (2,3)-coboundaries is contained in space of (2,3)-cocycles. The proof is finished.
\qed

\begin{defi}\label{def:1cohmologygroup}
The (2,3)-cohomology group  of a HLYA  $T$ with coefficients in $V$ is defined as the quotient space
$$H^{2}(T, V)\times H^{3}(T, V)\triangleq Z^{2}(T, V)\times Z^{3}(T, V)/B^{2}(T, V)\times B^{3}(T, V).$$
\end{defi}

\section{Infinitesimal Deformations}
\label{sec:3}

Let $T$ be a HLYA and $\nu: T\times T\to T$ and $\omega: T\times T\times T\to T$ be bilinear and trilinear maps. Consider a $\lambda$-parametrized family of bilinear maps and trilinear maps:
\begin{eqnarray*}
[x_1, x_2]_\lam&\triangleq& [x_1, x_2]+ \lambda\nu(x_1, x_2),\\
{[x_1, x_2, x_3]}_\lam&\triangleq& [x_1, x_2, x_3]+ \lambda\omega(x_1, x_2, x_3).
 \end{eqnarray*}

If $[\cdot,\cdot]_\lam$ and $[\cdot,\cdot,\cdot]_\lam$ endow $T$ with a HLYA structure which is denoted by $T_\lam$, then we say that $(\nu,\omega)$ generates a
$\lambda$-parameter infinitesimal deformation of HLYA $T$.

\begin{thm}\label{thm:deformation}
With the above notations, $(\nu,\omega)$ generates a $\lambda$-parameter infinitesimal deformation of a HLYA $T$ if and only if the following two conditions hold:

(i) $(\nu,\omega)$  defines a HLYA of deformation type on $T$;

(ii) $(\nu,\omega)$ is a (2,3)-cocycle of $T$ with coefficients in the adjoint representation.
\end{thm}

\pf
Assume $(\nu,\omega)$ generates a $\lambda$-parameter infinitesimal deformation of the HLYA $T$, then the
maps $[x_1, x_2]_\lam$ and ${[x_1, x_2, x_3]}_\lam$ defined above must satisfies conditions (HLY1)--(HLY6).
From these conditions, we will deduce that $(\nu,\omega)$ is a (2,3)-cocycle and $(\nu,\omega)$  defines a HLYA of deformation type on $T$.

From (HLY01), we have
\begin{eqnarray*}
&&\al([x_1, x_2]\dlam)-[\al (x_1),\al (x_2)]\dlam\\
&=&\al[x_1, x_2]-[\al (x_1),\al (x_2)]+ \lambda\{\al\circ\nu(x_1, x_2)-\nu(\al (x_1),\al (x_2))\}\\
&=&0.
\end{eqnarray*}
thus we get
\begin{eqnarray}\label{eq:dm001}
\al\circ\nu(x_1, x_2)=\nu(\al (x_1),\al (x_2)).
\end{eqnarray}

From (HLY02), we have
\begin{eqnarray*}
&&\al([x_1, x_2,x_3]\dlam)-[\al (x_1), \al (x_2),\al (x_3)]\dlam\\
&=&\al[x_1, x_2,x_3]-[\al (x_1),\al(x_2),\al(x_3)]\\
&&+\lambda\{\al\circ\nu(x_1, x_2,x_3)-\nu(\al (x_1),\al (x_2),\al(x_3))\}\\
&=&0.
\end{eqnarray*}
thus we obtain
\begin{eqnarray}\label{eq:dm002}
\al\circ\nu(x_1, x_2,x_3)=\nu(\al (x_1),\al (x_2),\al(x_3)).
\end{eqnarray}

From (HLY3), we have
\begin{eqnarray*}
&&[x_1, x_2,x_3]\dlam+c.p.+[[x_1, x_2]\dlam,\al(x_3)]\dlam+c.p.\\
&=&[x_1, x_2,x_3]+c.p.+[[x_1, x_2],\al(x_3)]+c.p.\\
&&+\lam\{\omega(x_1, x_2,x_3)+c.p.+\nu([x_1, x_2],\al(x_3))+c.p.+[\nu(x_1, x_2),\al(x_3)]+c.p.\}\\
&&+\lam^2\{\nu(\nu(x_1, x_2),\al(x_3))+c.p.\}\\
&=&0,
\end{eqnarray*}
thus we get
\begin{eqnarray}
\label{eq:dm03}\omega(x_1, x_2,x_3)+c.p.+\nu([x_1, x_2],\al(x_3))+c.p.+[\nu(x_1, x_2),\al(x_3)]+c.p.=0,\\
\label{eq:dm03'}\nu(\nu(x_1, x_2),\al(x_3))+c.p.=0.
\end{eqnarray}

From (HLY4), we have
\begin{eqnarray*}
&&[[x_1, x_2]\dlam, \al(x_3), \al(y_1)]\dlam+c.p.\\
&=&[[x_1, x_2], \al(x_3), \al(y_1)]+c.p\\
&&\lam\{\omega([x_1, x_2],\al(x_3),\al(y_1))+c.p.+[\nu(x_1, x_2),\al(x_3),\al(y_1)]+c.p.\}\\
&&+\lam^2\{\omega(\nu(x_1, x_2),\al(x_3),\al(y_1))+c.p.\}\\
&=&0,
\end{eqnarray*}
thus we get
\begin{eqnarray}
\label{eq:dm04}\omega([x_1, x_2],\al(x_3),\al(y_1))+c.p.+[\nu(x_1, x_2),\al(x_3),\al(y_1)]+c.p.=0,\\
\label{eq:dm04'}\omega(\nu(x_1, x_2),\al(x_3),\al(y_1))+c.p.=0.
\end{eqnarray}

From (HLY5), we have
\begin{eqnarray*}
[\al(x_1), \al(x_2), [y_1, y_2]\dlam]\dlam =[[ x_1, x_2, y_1]\dlam ,  \al^2(y_2)]\dlam  + [\al^2(y_1), [ x_1, x_2, y_2]\dlam]\dlam,
\end{eqnarray*}
the left hand side is equal to
\begin{eqnarray*}
&&[\al(x_1), \al(x_2), [y_1, y_2]+\lam\nu(y_1, y_2)]\dlam\\
&=&[\al(x_1), \al(x_2), [y_1, y_2]]\\
&&+\lam\{\omega(\al(x_1), \al(x_2), [y_1, y_2])+[\al(x_1), \al(x_2),\nu(y_1, y_2)]\}\\
&&+\lam^2\omega(\al(x_1), \al(x_2),\nu(y_1, y_2)),
\end{eqnarray*}
and the right hand side is equal to
\begin{eqnarray*}
&&[[ x_1, x_2, y_1]\dlam , \al^2(y_2)]\dlam  + [\al^2(y_1), [ x_1, x_2, y_2]\dlam]\dlam\\
&=&[[ x_1, x_2, y_1] , \al^2(y_2)] + [\al^2(y_1), [ x_1, x_2, y_2]]\\
&&+\lam\{[\omega(x_1, x_2, y_1), \al^2(y_2)]+\nu([x_1, x_2,y_1], \al^2(y_2))\\
&&+[\al^2(y_1), \omega(x_1, x_2, y_2]+\nu(y_1, [x_1, x_2,y_2])\}\\
&&+\lam^2\{\nu(\omega( x_1, x_2, y_1), \al^2(y_2)) + \nu(\al^2(y_1),\omega( x_1, x_2, y_2))\},
\end{eqnarray*}
then we obtain
\begin{eqnarray}
\notag&&\omega(\al(x_1), \al(x_2), [y_1, y_2])+[\al(x_1), \al(x_2),\nu(y_1, y_2)]\\
\notag&=&[\omega(\al(x_1), \al(x_2), y_1), \al^2(y_2)]+\nu([x_1, x_2,y_1], \al^2(y_2))\\
\label{eq:dm05}&&+[\al^2(y_1), \omega(x_1, x_2, y_2)]+\nu(\al^2(y_1), [x_1, x_2,y_2]),
\end{eqnarray}
and
\begin{eqnarray}
\label{eq:dm05'}\omega(\al(x_1), \al(x_2),\nu(y_1, y_2))&=&\nu(\omega( x_1, x_2, y_1), \al^2(y_2)) + \nu(\al^2(y_1), \omega( x_1, x_2, y_2)).
\end{eqnarray}

From (HLY6), we have
\begin{eqnarray*}
&&[\al^2(x_1), \al^2(x_2), [y_1, y_2, y_3]\dlam]\dlam \\
&=&[[ x_1, x_2, y_1]\dlam , \al^2(y_2), \al^2(y_3)]\dlam  + [\al^2(y_1),[ x_1, x_2, y_2]\dlam , \al^2(y_3)]\dlam \\
&& + [\al^2(y_1), \al^2(y_2),[ x_1, x_2, y_3]\dlam]\dlam,
\end{eqnarray*}
the left hand side is equal to
\begin{eqnarray*}
&&[\al^2(x_1), \al^2(x_2),[y_1, y_2, y_3]+\lam\omega(y_1, y_2, y_3)]\dlam\\
&=&[\al^2(x_1), \al^2(x_2), [y_1, y_2, y_3]]+\lam\omega(\al^2(x_1), \al^2(x_2), [y_1, y_2, y_3])\\
&&+[\al^2(x_1), \al^2(x_2),\lam\omega(y_1, y_2, y_3)]+\lam\omega(\al^2(x_1), \al^2(x_2),\lam\omega(y_1, y_2, y_3))\\
&=&[\al^2(x_1), \al^2(x_2), [y_1, y_2, y_3]]\\
&&+\lam\{\omega(\al^2(x_1), \al^2(x_2), [y_1, y_2, y_3])+[\al^2(x_1), \al^2(x_2),\omega(y_1, y_2, y_3)]\}\\
&&+\lam^2\omega(\al^2(x_1), \al^2(x_2),\omega(y_1, y_2, y_3)),
\end{eqnarray*}
and the right hand side is equal to
\begin{eqnarray*}
&&[[x_1, x_2, y_1]+\lam\omega(x_1, x_2, y_1),\al^2(y_2), \al^2(y_3)]\dlam \\
&& + [\al^2(y_1),  [ x_1, x_2, y_2]+\lam\omega(x_1, x_2, y_2), \al^2(y_3)]\dlam \\
&&+ [\al^2(y_1), \al^2(y_2),  [ x_1, x_2, y_3]+\lam\omega(x_1, x_2, y_3)]\dlam\\
&=&[[ x_1, x_2, y_1] , \al^2(y_2), \al^2(y_3)]  + [y_1, [ x_1, x_2, y_2] , y_3] + [\al^2(y_1), \al^2(y_2), [ x_1, x_2, y_3]]\\
&&+\lam\{\omega([x_1, x_2, y_1],\al^2(y_2), \al^2(y_3)) +[\omega(x_1, x_2, y_1), \al^2(y_2), \al^2(y_3)]\\
&&\qquad+\omega(\al^2(y_1),[ x_1, x_2, y_2],\al^2(y_3))+[\al^2(y_1),\omega(x_1, x_2, y_2),\al^2(y_3)]\\
&&\qquad+\omega(\al^2(y_1), \al^2(y_2), [ x_1, x_2, y_3])+[\al^2(y_1), \al^2(y_2),\omega(x_1, x_2, y_3)] \}\\
&&+\lam^2\{\omega(\omega(x_1, x_2, y_1), \al^2(y_2), \al^2(y_3)) + \omega(y_1, \omega( x_1, x_2, y_2), y_3) \\
&&+ \omega(\al^2(y_1), \al^2(y_2), \omega( x_1, x_2, y_3))\},
\end{eqnarray*}
then we get
\begin{eqnarray}
\nonumber &&\omega(\al^2(x_1), \al^2(x_2),[y_1, y_2, y_3])+[\al^2(x_1), \al^2(x_2),\omega(y_1, y_2, y_3)]\\
\nonumber &=&\omega([x_1, x_2, y_1],\al^2(y_2), \al^2(y_3))+\omega(\al^2(y_1),[ x_1, x_2, y_2], \al^2(y_3)) \\
\nonumber &&+\omega(\al^2(y_1), \al^2(y_2),[ x_1, x_2, y_3])+[\omega(x_1, x_2, y_1),\al^2(y_2), \al^2(y_3)]\\
\label{eq:dm06} &&+[\al^2(y_1),\omega(x_1, x_2, y_2),\al^2(y_3)]+[\al^2(y_1), \al^2(y_2),\omega(x_1, x_2, y_3)],
\end{eqnarray}
and
\begin{eqnarray}
\nonumber &&\omega(\al^2(x_1), \al^2(x_2),\omega(y_1, y_2, y_3))\\
\nonumber&=&\omega(\omega( x_1, x_2, y_1),\al^2(y_2), \al^2(y_3)) + \omega(y_1, \omega( x_1, x_2, y_2), y_3)\\
\label{eq:dm06'}&&+ \omega(\al^2(y_1), \al^2(y_2), \omega( x_1, x_2, y_3)).
\end{eqnarray}
Therefore by \eqref{eq:dm001}, \eqref{eq:dm002}, \eqref{eq:dm03'}, \eqref{eq:dm04'}, \eqref{eq:dm05'} and \eqref{eq:dm06'}, $(\nu,\omega)$ defines a HLYA of deformation type on $T$.
Furthermore, by \eqref{eq:dm001}, \eqref{eq:dm002}, \eqref{eq:dm03}, \eqref{eq:dm04}, \eqref{eq:dm05} and \eqref{eq:dm06}, we obtain that $(\nu,\omega)$ is a (2,3)-cocycle of $T$ with coefficients in the adjoint representation.
\qed
\medskip

A deformation is said to be {\bf trivial} if there exists a linear map $N: T\to  T$
such that for $\varphi_\lam = \id + \lambda N$: $T_\lam \to  T$ there hold
\begin{eqnarray}\label{eq:Nijenhuis00}
\varphi_\lam [x_1,x_2]\dlam=[\varphi_\lam x_1,\varphi_\lam x_2]\quad\mbox{and}\quad
\varphi_\lam [x_1,x_2,x_3]\dlam=[\varphi_\lam x_1,\varphi_\lam x_2,\varphi_\lam x_3].
\end{eqnarray}
It follows from \eqref{eq:Nijenhuis00} that $N$ must satisfy the following condition
\begin{eqnarray}
\label{eq:Nijenhuis03}N[Nx_1,x_2]+N[x_1,Nx_2]-N^2[x_1,x_2]=[Nx_1,Nx_2];
\end{eqnarray}
and
\begin{eqnarray}
\label{eq:Nijenhuis4}\nonumber &&N[Nx_1,x_2,x_3]+N[x_1,Nx_2,x_3]+N[x_1,x_2,Nx_3]-N^2[x_1,x_2,x_3]\\
&=&[Nx_1,Nx_2,x_3]+[Nx_1,x_2,Nx_3]+[x_1,Nx_2,Nx_3].
\end{eqnarray}

\emptycomment{
\begin{eqnarray}\label{eq:Nijenhuis4'}
\nonumber N^2[x_1,x_2,x_3]&=&N[Nx_1,x_2,x_3]+N[x_1,Nx_2,x_3]+N[x_1,x_2,Nx_3]\\
&&-([Nx_1,Nx_2,x_3]+[Nx_1,x_2,Nx_3]+[x_1,Nx_2,Nx_3]).
\end{eqnarray}}

\begin{defi}
A linear operator $N: T\to  T$ is called a Nijenhuis operator of a HLYA $T$
if \eqref{eq:Nijenhuis03} and \eqref{eq:Nijenhuis4} hold.
\end{defi}

An important property of Nijenhuis operator is that it gives trivial deformation.

\begin{thm}Let $N$ be a Nijenhuis operator for $T$. Then a deformation of $T$ can be obtained by putting
\begin{eqnarray}
\label{eq:bracket01}\nu(x_1,x_2)&=&[Nx_1,x_2]+[x_1,Nx_2]-N[x_1,x_2],\\
\label{eq:bracket02}\omega(x_1,x_2,x_3)&=&[Nx_1,x_2,x_3]+[x_1,Nx_2,x_3]+[x_1,x_2,Nx_3]-N[x_1,x_2,x_3].
\end{eqnarray}
Furthermore, this deformation is a trivial one.
\end{thm}

\section{Abelian Extensions} 
\label{sec:4}

In this section, we study abelian extensions of HLYAs.
It is showed that abelian extensions are classified by the (2,3)-cohomology group.
We will built a bijection map from the set of equivalent classes of abelian extensions $\Ext( T,\frkh)$ and
$H^2(T, V)\times H^3(T, V)$.

\begin{defi}
 Let $( T,[\cdot,\cdot], [\cdot,\cdot,\cdot],\al)$, $(\frkh,[\cdot,\cdot]_\frkh, [\cdot,\cdot,\cdot]_\frkh,\be)$ and
 $(\hat{T}, [\cdot,\cdot]_{\hat{T}}, [\cdot,\cdot,\cdot]_{\hat{T}},\hal)$ be HLYAs,
$i:\frkh\to\hat{T},~~p:\hat{T}\to T$
be homomorphisms. If  the following diagram commutes and the horizontal two lines are
short exact sequence (i.e.$\mathrm{Im}(i)=\mathrm{Ker}(p)$,$\mathrm{Ker}(i)=0$ and $\mathrm{Im}(p)= T$),
\begin{equation}\label{diagram:exact}
 \xymatrix{
   0  \ar[r]^{} & \frkh \ar[d]_{\be} \ar[r]^{i} & \hat{T} \ar[d]_{\hal} \ar[r]^{p} &  T \ar[d]_{\al} \ar[r]^{} & 0 \\
   0 \ar[r]^{} & \frkh \ar[r]^{i} & \tilde{T} \ar[r]^{p} &  T \ar[r]^{} & 0.
  }
\end{equation}
then we call $\hat{T}$  an extension of $T$ through
$\frkh$, and denote it by $\E_{\hat{T}}$.
It is called an abelian extension if $\frkh$ is abelian ideal of $\hat{T}$, i.e. $[u,v]_{\hat{T}}=0$ and
$[u,v,\cdot]_{\hat{T}}=[u,\cdot,v]_{\hat{T}}=[\cdot,u,v]_{\hat{T}}=0$,
for all $u,v\in \frkh$.
\end{defi}

From the left square in the commutative diagram we deduce that if we choose element $u\in V$, then
\begin{equation}
\hal\circ i(u)=i\circ \al_v(u).
\end{equation}
Since $i$ is a injective map,  we can identify $V$ with its image in $\hat{T}$, thus we have
\begin{equation}
\hal(u)=\hal|_{V}(u)=\be(u).
\end{equation}

\emptycomment{
This mean that $\hat{T}$ is a HLYA containing $V$ as an ideal such that $\hat{T}/V\simeq T$.
thus $\hat{T}\simeq T\oplus V$ as vector space.
$$\hal\circ i(u)=i\circ \al_v(u).$$
$$\hal(x)= p(\hal(x)+\hal(u))=p\circ\hal(x+u)=\al\circ p(x+u)=\al(x).$$
$$\hal|_{V}(u)=\al_v(u),\quad \hal|_{T}(x)=\al(x).$$
$$\hal(u)=\al_v(u),\quad \hal(x)=\al(x).$$
}

A section $\sigma: T\to\hat{T}$ of $p:\hat{T}\to T$
consists of linear maps $\sigma: T\to\hat{T}$ such that
\begin{equation}
p\circ\sigma=\id_{T}\,\ \mbox{and}\,\ \hal\circ\sigma=\sigma\circ\al.
\end{equation}

\begin{defi}
 Two extensions of HLYAs
 $\E_{\hat{T}}:0\to\frkh\stackrel{i}{\to}\hat{T}\stackrel{p}{\to} T\to0$
 and $\E_{\tilde{T}}:0\to\frkh\stackrel{j}{\to}\tilde{T}\stackrel{q}{\to} T\to0$ are called equivalent,
 if there exists a HLYA homomorphism $F:\hat{T}\to\tilde{T}$  such that the following diagram commutes
\begin{equation}\label{diagram:equivalent}
\xymatrix{
   0  \ar[r]^{} & \frkh \ar[d]_{\id} \ar[r]^{i} & \hat{T} \ar[d]_{F} \ar[r]^{p} &  T \ar[d]_{\id} \ar[r]^{} & 0 \\
   0 \ar[r]^{} & \frkh \ar[r]^{j} & \tilde{T} \ar[r]^{q} &  T \ar[r]^{} & 0.
   }
\end{equation}
The set of equivalent classes of extensions of $T$ by $\h$ is denoted by $\Ext( T,\frkh)$.
\end{defi}

Let $\hat{T}$ be an abelian extension of $T$ by $\frkh$.  Define maps $\rho$ from $T$ to $\End(\frkh)$ and
$\DD, \theta$ from $T\times  T$ to $\End(\frkh)$ by
\begin{eqnarray}
\label{eq:rep01}\rho(x_1)(u)&\triangleq&[\sigma(x_1),u]_{\hat{T}},\\ 
\label{eq:rep02}\DD(x_1,x_2)(u)&\triangleq&[\sigma(x_1),\sigma(x_2),u]_{\hat{T}},\\
\label{eq:rep03}\theta(x_1,x_2)(u)&\triangleq&[u,\sigma(x_1),\sigma(x_2)]_{\hat{T}}.
\end{eqnarray}

\begin{lem}\label{lem:11}
With the above notations, $(\rho,\DD,\theta)$ is a representation of $T$ on $V$ and does not depend on the choice of the section $\sigma$.
Moreover,  equivalent abelian extensions give the same representation.
\end{lem}

\pf
First, the fact that $\rho,  \DD, \theta$ are independent of
the choice of $\sigma$ is easy to check. For details, see \cite{Zhang1}.

Second, we will show that $(\rho,\DD,\theta)$ is a representation of $T$ on $V$.

By the equality
\begin{eqnarray*}
\rho(\al(x_1))\circ \be(u)&=&{[\sigma(\al(x_1)),\be(u)]_{\hat{T}}}\\
&=&{[\hal\circ\sigma(x_1),\hal(u)]_{\hat{T}}}\\
&=&\hal([\sigma(x_1),u]_{\hat{T}})\\
&=&\be([\sigma(x_1),u]_{\hat{T}})=\be\circ\rho(x_1)(u),
\end{eqnarray*}
we obtain (HR01):
\begin{eqnarray}
\rho(\al(x_1))\circ \be=\be\circ\rho(x_1).
\end{eqnarray}

By the equality
\begin{eqnarray*}
\DD(\al(x_1),\al(x_2))\circ \be(u)&=&{[\sigma(\al(x_1)),\sigma(\al(x_2)),\be(u)]_{\hat{T}}}\\
&=&{[\hal\circ\sigma(x_1),\hal\circ\sigma(x_2),\hal(u)]_{\hat{T}}}\\
&=&\hal([\sigma(x_1),\sigma(x_2),u]_{\hat{T}})\\
&=&\be([\sigma(x_1),\sigma(x_2),u]_{\hat{T}})=\be\circ\DD(x_1,x_2)(u),
\end{eqnarray*}
we obtain (HR02):
\begin{eqnarray}
\DD(\al(x_1),\al(x_2))\circ \be(u)=\be\circ\DD(x_1,x_2)(u).
\end{eqnarray}

By the equality
\begin{eqnarray*}
&&[\si (x_1),\si (x_2),u]_{\hat{T}}+[\si (x_2), u,\si (x_1)]_{\hat{T}}+[u, \si (x_1),\si (x_2)]_{\hat{T}}\\
&&+[[\si (x_1),\si (x_2)]_{\hat{T}},\be(u)]_{\hat{T}}+[[\si (x_2), u]_{\hat{T}},\hal\circ\si (x_1)]_{\hat{T}}+[[u, \si (x_1)],\hal\circ\si (x_2)]_{\hat{T}}=0,
\end{eqnarray*}
we obtain (HR31):
\begin{eqnarray}
\DD(x_1,x_2)-\theta(x_2,x_1)+\theta(x_1,x_2)+\rho([x_1,x_2])\circ\be-\rho(\al(x_1))\rho (x_2)+\rho (\al(x_2))\rho (x_1)=0.
\end{eqnarray}

By the equality
\begin{eqnarray*}
&&[[\si (x_1),\si (x_2)]_{\hat{T}},\hal\circ\si (x_3),\be(u)]_{\hat{T}}+[[\si (x_2), \si (x_3)]_{\hat{T}},\hal\circ\si (x_1),\be(u)]_{\hat{T}}\\
&&+[[\si (x_3), \si (x_1)]_{\hat{T}},\hal\circ\si (x_2),\be(u)]_{\hat{T}}=0,
\end{eqnarray*}
we have (HR41):
\begin{eqnarray}
\DD([x_1,x_2],\al(x_3))+\DD([x_2,x_3],\al(x_1))+\DD([x_3,x_1],\al(x_2))=0.
\end{eqnarray}

By the equality
\begin{eqnarray*}
&&[[\si (x_1),\si (x_2)]_{\hat{T}},\be(u),\hal\circ\si (y_1)]_{\hat{T}}+[[\si (x_2), u]_{\hat{T}},\hal\circ\si (x_1),\hal\circ\si (y_1)]_{\hat{T}}\\
&&+[[u, \si (x_1)]_{\hat{T}},\hal\circ\si (x_2),\hal\circ\si (y_1)]_{\hat{T}}=0,
\end{eqnarray*}
we have (HR42):
\begin{eqnarray}
\theta([x_1,x_2],\al(y_1))\circ \be=\theta(\al(x_1),\al(y_1))\rho(x_2)-\theta(\al(x_2),\al(y_1))\rho(x_1).
\end{eqnarray}

By the equality
\begin{eqnarray*}
[\hal\circ\si (x_1),\hal\circ\si (x_2), [u, \si (y_2)]_{\hat{T}}]_{\hat{T}} &=&
[[\si (x_1),\si (x_2),u]_{\hat{T}},\hal^2\circ\si  (y_2)]_{\hat{T}} \\
&&+ [\be^2(u), [\si (x_1),\si (x_2),\si (y_2)]_{\hat{T}}]_{\hat{T}},
\end{eqnarray*}
we have (HR51):
\begin{eqnarray}
&&\DD(\al(x_1),\al(x_2))\rho(y_2)=\rho(\al^2(y_2))\DD(x_1,x_2)+\rho([x_1, x_2, y_2])\circ\bebe.
\end{eqnarray}

By the equality
\begin{eqnarray*}
[\hal\circ\si (x_1), \be(u), [\si (y_1), \si (y_2)]_{\hat{T}}]_{\hat{T}} &=& [[\si (x_1),u,\si (y_1)]_{\hat{T}},\hal^2\circ\si  (y_2)]_{\hat{T}}\\
&& + [\hal^2\circ\si (y_1), [\si (x_1),u,\si (y_2)]_{\hat{T}}]_{\hat{T}},
\end{eqnarray*}
we have (HR52):
\begin{eqnarray}
\theta(\al(x_1),[y_1, y_2])\circ \be=\rho(\al^2(y_1))\theta(x_1, y_2)-\rho(\al^2(y_2))\theta(x_1,y_1).
\end{eqnarray}

By the equality
\begin{eqnarray*}
&&[\hal^2\circ\si (x_1),\hal^2\circ\si (x_2), [u,\si (y_1),\si (y_2)]_{\hat{T}}]_{\hat{T}} \\
&=& [[\si (x_1),\si (x_2),u]_{\hat{T}},\hal^2\circ\si (y_1),\hal^2\circ\si (y_2)]_{\hat{T}} + [\be^2(u), [\si (x_1),\si (x_2),\si (y_1)]_{\hat{T}},\hal^2\circ\si (y_2)]_{\hat{T}}\\
&&+ [\be^2(u),\hal^2\circ\si (y_1), [\si (x_1),\si (x_2),\si (y_2)]_{\hat{T}}]_{\hat{T}},
\end{eqnarray*}
we have (HR61):
\begin{eqnarray}
\nonumber\DD(\al^2(x_1),\al^2(x_2))\theta(y_1,y_2)
&=&\theta(\al^2(y_1),\al^2(y_2))\DD(x_1,x_2)+\theta([x_1,x_2,y_1],\al^2(y_2))\circ\be^2\\
&&+\theta(\al^2(y_1),[x_1,x_2,y_2])\circ\be^2.
\end{eqnarray}

By the equality
\begin{eqnarray*}
&&[\be^2(u),\hal^2\circ\si (x_1), [\si (y_1),\si (y_2),\si (y_3)]_{\hat{T}}]_{\hat{T}} \\
&=& [[u,\si (x_1),\si (y_1)]_{\hat{T}},\hal^2\circ\si (y_2),\hal^2\circ\si (y_3)]_{\hat{T}}
+ [\hal^2\circ\si (y_1), [u,\si (x_1),\si (y_2)]_{\hat{T}},\hal^2\circ\si (y_3)]_{\hat{T}}\\
&&+ [\hal^2\circ\si (y_1),\hal^2\circ\si (y_2), [u,\si (x_1),\si (y_3)]_{\hat{T}}]_{\hat{T}},
\end{eqnarray*}
we have (HR62):
\begin{eqnarray}
\nonumber\theta(\al^2(x_1),[y_1, y_2, y_3])\circ\be^2&=& \theta(\al^2(y_2), \al^2(y_3))\theta(x_1,y_1) - \theta(\al^2(y_1), \al^2(y_3))\theta(x_1,y_2)\\
&&+ \DD (\al^2(y_1), \al^2(y_2))\theta(x_1,y_3).
\end{eqnarray}
Therefore we see that $(\rho,  \DD, \theta)$ is a representation of $T$ on $V$.

At last, suppose that $\E_{\hat{T}}$ and $\E_{\tilde{T}}$ are equivalent abelian extensions, and $F:\hat{T}\to\tilde{T}$ is the HLYA homomorphism satisfying $F\circ i=j$, $q\circ F=p$.
Choose linear sections $\sigma$ and $\sigma'$ of $p$ and $q$, we get $qF\sigma(x_i)=p\sigma(x_i)=x_i=q\sigma'(x_i)$,
then $F\sigma(x_i)-\sigma'(x_i)\in \Ker (q)\cong\h$. Thus, we have
$$
[u,\sigma(x_1),\sigma(x_2)]_{\hat{T}}=[u,F\sigma(x_1),F\sigma(x_2)]_{\tilde{T}}=[u,\sigma'(x_1),\sigma'(x_2)]_{\tilde{T}}.
$$
Therefore, equivalent abelian extensions give the same $\theta$.
Similarly, one can prove that equivalent abelian extensions give the same $\DD$ and $\rho$.
The proof is finished.
\qed

\medskip

Let $\sigma: T\to\hat{T}$  be a section of abelian extension. Define the following maps:
\begin{eqnarray}
\label{eq:coc01}\nu(x_1,x_2)&\triangleq&[\sigma(x_1),\sigma(x_2)]_{\hat{T}}-\sigma([x_1,x_2]),\\
\label{eq:coc02}\omega(x_1,x_2,x_3)&\triangleq&[\sigma(x_1),\sigma(x_2),\sigma(x_3)]_{\hat{T}}-\sigma([x_1,x_2,x_3]),
\end{eqnarray}

\begin{lem}\label{lem:22}
Let $0\to\frkh{\to}\hat{T}{\to} T\to 0$ be an abelian extension of $T$ by $\h$.
Then $(\nu,\omega)$ defined by \eqref{eq:coc01} and \eqref{eq:coc02} is a $(2,3)$-cocycle of $T$ with coefficients in $\frkh$.
\end{lem}

\pf First, we claim that the image of $\nu$ is contained in $V$, that is to say, $p\circ \nu(x_1,x_2)=0$.
In fact, since $p$ is an algebraic homomorphism, we have
$$p\circ \nu(x_1,x_2)=[p\circ\sigma(x_1),p\circ\sigma(x_2)]_{\hat{T}}-p\circ\sigma([x_1,x_2])=0$$

Next, one check that $\nu$ and $\omega$ defined above satisfies (CC01) and (CC02). For example
\begin{eqnarray*}
&&\nu(\al(x_1),\al(x_2))\\
&=&[\sigma(\al(x_1)),\sigma(\al(x_2))]_{\hat{T}}-\sigma([\al(x_1),\al(x_2)])\\
&=&{[\sigma\circ\al(x_1),\sigma\circ\al(x_2)]_{\hat{T}}-\sigma\circ\al([x_1,x_2])}\\
&=&{([\hal\circ\sigma(x_1),\hal\circ\sigma(x_2)]_{\hat{T}})-\hal\circ\sigma([x_1,x_2])}\\
&=&\hal([\sigma(x_1),\sigma(x_2)]_{\hat{T}})-\sigma([x_1,x_2]))\\
&=&\be([\sigma(x_1),\sigma(x_2)]_{\hat{T}}-\sigma([x_1,x_2]))\\
&=&\be(\nu(x_1,x_2)).
\end{eqnarray*}


Finally, we verify that $\nu$ and $\omega$ satisfies (CC1)--(CC4).

By the equality
\begin{eqnarray*}
&&[\si x_1,\si x_2,\si x_3]_{\hat{T}}+c.p.+[[\si x_1,\si x_2]_{\hat{T}},\hal(\si (x_3))]_{\hat{T}}+c.p.=0,
\end{eqnarray*}
we obtain that
\begin{eqnarray*}
&&\{\omega([x_1,x_2, x_3)+\si[x_1,x_2,x_3]_{\hat{T}}\}+c.p.\\
&&+\{[\nu(x_1,x_2),\si (\al(x_3))]_{\hat{T}}+\nu([x_1,x_2],\al(x_3))+\si([[x_1,x_2], \al(x_3)])\}+c.p.=0.
\end{eqnarray*}
Thus we have (CC1):
\begin{eqnarray}
&&\omega(x_1,x_2, x_3)+c.p.-\rho(\al(x_3))\nu(x_1,x_2)-c.p.+\nu([x_1,x_2],\al(x_3))+c.p.=0.
\end{eqnarray}

By the equality
\begin{eqnarray*}
&&[[\si x_1,\si x_2]_{\hat{T}},\hal(\si x_3),\hal(\si y_1)]_{\hat{T}}+c.p.=0,
\end{eqnarray*}
we get
\begin{eqnarray*}
&&\{[\nu(x_1,x_2),\si \al(x_3),\si \al(y_1)]_{\hat{T}}+\omega([x_1,x_2], \al(x_3), \al(y_1))+\si[[x_1,x_2], \al(x_3), \al(y_1)]\}+c.p.=0.
\end{eqnarray*}
Thus we have (CC2):
\begin{eqnarray}
&&\theta(\al(x_3), \al(y_1))\nu(x_1,x_2)+c.p.+\omega([x_1,x_2],\al(x_3), \al(y_1))+c.p.=0.
\end{eqnarray}

By the equality
\begin{eqnarray*}
&&[\hal(\si x_1),\hal(\si x_2), [\si y_1,\si y_2]_{\hat{T}}]_{\hat{T}} \\
&=& [[\si x_1,\si x_2,\si y_1]_{\hat{T}},\hal^2(\si y_2)]_{\hat{T}} + [\hal^2(\si y_1), [\si x_1,\si x_2,\si y_2]_{\hat{T}}]_{\hat{T}},
\end{eqnarray*}
we obtain that the left hand side is equal to
\begin{eqnarray*}
&&[\si\al(x_1),\si \al(x_2), [\si y_1,\si y_2,]_{\hat{T}}]_{\hat{T}}\\
&=&[\si\al(x_1),\si \al(x_2), \nu(y_1,y_2)+\sigma([y_1,y_2]_ T)]_{\hat{T}}\\
&=&\DD(\al(x_1), \al(x_2))\nu(y_1,y_2)+[\si\al(x_1),\si \al(x_2),\sigma([y_1,y_2])]_{\hat{T}}\\
&=&\DD(\al(x_1), \al(x_2))\nu(y_1,y_2)+\omega(\al(x_1), \al(x_2),[y_1,y_2])+\sigma([\al(x_1), \al(x_2),[y_1,y_2]]).
\end{eqnarray*}
Similarily, the right hand side is equal to
\begin{eqnarray*}
&&[[\si x_1,\si x_2,\si y_1]_{\hat{T}},\si \al^2(y_2)]_{\hat{T}} + [\si \al^2(y_1), [\si x_1,\si x_2,\si y_2]_{\hat{T}}]_{\hat{T}}\\
&=&[\omega(x_1,x_2,y_1)+\si[x_1,x_2,y_1],\si y_2]_{\hat{T}}+\nu([x_1,x_2,y_1],y_2)+\sigma[[x_1,x_2,y_1],y_2]\\\
&=&-\rho(\al^2(y_2))\omega(x_1,x_2,y_1)+\nu([ x_1, x_2, y_1], \al^2(y_2))+\sigma([[ x_1, x_2, y_1],\al^2(y_2)])\\
&&\rho(\al^2(y_1))\omega(x_1,x_2,y_2)+\nu(\al^2(y_1), [ x_1, x_2, y_2])+\sigma([\al^2(y_1), [ x_1, x_2, y_2]]).
\end{eqnarray*}
Thus we have (CC3):
\begin{eqnarray}
\nonumber&& \omega(\al(x_1), \al(x_2),[y_1, y_2])+\DD(\al(x_1), \al(x_2))\omega(y_1, y_2)\\
\nonumber&=&\nu([x_1, x_2, y_1], \al^2(y_2)) + \nu(\al^2(y_1),[x_1, x_2,y_2]))\\
&&+\rho(\al^2(y_1))\omega(x_1,x_2,y_2)-\rho(\al^2(y_2))\omega(x_1,x_2,y_1).
\end{eqnarray}

By the equality
\begin{eqnarray*}
&&[\hal^2(\si x_1),\hal^2(\si x_2), [\si y_1,\si y_2,\si y_3]_{\hat{T}}]_{\hat{T}} \\
&=& [[\si x_1,\si x_2,\si y_1]_{\hat{T}},\hal^2(\si y_2),\hal^2(\si y_3)]_{\hat{T}} + [\hal^2(\si y_1), [\si x_1,\si x_2,\si y_2]_{\hat{T}},\hal^2(\si y_3)]_{\hat{T}}\\
&& + [\hal^2(\si y_1),\hal^2(\si y_2), [\si x_1,\si x_2,\si y_3]_{\hat{T}}]_{\hat{T}},
\end{eqnarray*}
we have that the left hand side is equal to
\begin{eqnarray*}
&&[\si \al^2(x_1),\si\al^2(x_2), [\si y_1,\si y_2,\si y_3]_{\hat{T}}]_{\hat{T}}\\
&=&[\si \al^2(x_1),\si\al^2(x_2), \omega(y_1,y_2,y_3)+\sigma([y_1,y_2,y_3]g)]_{\hat{T}}\\
&=&\DD(\al^2(x_1), \al^2(x_2))\omega(y_1,y_2,y_3)+[\si x_1,\si x_2,\sigma([y_1,y_2,y_3])]_{\hat{T}}\\
&=&\DD(\al^2(x_1), \al^2(x_2))\omega(y_1,y_2,y_3)+\omega(\al^2(x_1), \al^2(x_2),[y_1,y_2,y_3])\\
&&+\sigma([\al^2(x_1),\al^2(x_2),[y_1,y_2,y_3]]).
\end{eqnarray*}
Similarily, the right hand side is equal to
\begin{eqnarray*}
&& \theta(y_2, y_3)\omega(x_1,x_2,y_1)+\omega([ x_1, x_2, y_1], y_2, y_3)+\sigma([[ x_1, x_2, y_1], y_2, y_3])\\
&&-\theta(y_1, y_3)\omega(x_1,x_2,y_2)+\omega(y_1, [ x_1, x_2, y_2], y_3)+\sigma([y_1, [ x_1, x_2, y_2], y_3])\\
&&+\DD(y_1, y_2) \omega(x_1,x_2,y_3)+ \omega(y_1, y_2, [ x_1, x_2, y_3])+\sigma([\al^2(y_1),\al^2(y_2), [ x_1, x_2, y_3]]).
\end{eqnarray*}
Thus we have (CC4):
\begin{eqnarray}
\nonumber&& \omega(\al^2(x_1), \al^2(x_2),[y_1, y_2, y_3])+\DD(\al^2(x_1), \al^2(x_2))\omega(y_1, y_2, y_3)\\
\nonumber&=&\omega([x_1, x_2, y_1], \al^2(y_2), \al^2(y_3)) + \omega(\al^2(y_1),[x_1, x_2,y_2], \al^2(y_3)) \\
\nonumber&&+ \omega(\al^2(y_1),\al^2(y_2), [x_1, x_2,y_3])+\theta(\al^2(y_2), \al^2(y_3))\omega(x_1,x_2,y_1)\\
&&- \theta(\al^2(y_1),\al^2(y_3))\omega(x_1,x_2,y_2) + \DD(\al^2(y_1),\al^2(y_2))\omega(x_1,x_2,y_3).
\end{eqnarray}
Therefore we get all the (2,3)-cocycle conditions in Definition \ref{def:1coc}.
\qed
\medskip

From Lemma \ref{lem:11} and Lemma \ref{lem:22}, we have obtained that abelian extensions of HLYA $T$ through $\frkh$ give rise to
a representation of $T$ on $V$ and a $(2,3)$-cocycle of $T$ with coefficients in $\frkh$.
Conversely, given a a representation and a $(2,3)$-cocycle,
we can obtain a HLYA structure on the space $T\oplus\frkh$.

\begin{lem}\label{lem:33}
Let $T$ be a HLYA,  $(\rho,\DD,\theta)$ is a representation of $T$ on  $V$ and $(\nu,\omega)$ is a (2,3)-cocycle of $T$ with coefficients in $\frkh$.
Then $T\oplus V$ is a HLYA under the following bilinear and trilinear maps:
\begin{eqnarray}
\label{eq:newbracket00}{(\al+\be)(x_1 + u_1)}&\triangleq&\al(x_1)+\be(u_1),\\
\label{eq:newbracket01}{[x_1 + u_1, x_2 + u_2]}_\nu&\triangleq&[x_1, x_2] + \nu(x_1, x_2)+\rho(x_1)(u_2) -\rho(x_2)(u_1),\\
\notag{[x_1 + u_1, x_2 + u_2, x_3 + u_3]}_\omega&\triangleq&[x_1, x_2, x_3] + \omega(x_1, x_2, x_3)+ \DD(x_1, x_2)(u_3)\\
\label{eq:newbracket02}&&-\theta(x_1, x_3)(u_2) + \theta( x_2, x_3)(u_1),
\end{eqnarray}
This kind of  HLYAs is denoted by $E_{(\nu,\omega)}=T\oplus_{(\nu,\omega)}\h$.
\end{lem}

\pf
We will verify that conditions (HLY01)--(HLY02) and (HLY1)--(HLY6) hold for maps defined on $T\oplus V$ by \eqref{eq:newbracket00},  \eqref{eq:newbracket01} and \eqref{eq:newbracket02}.

Now condition (HLY01) becomes
\begin{eqnarray}\label{semi01}
(\al+\be)([x_1 + u_1, x_2 + u_2])&=&[(\al+\be)(x_1 + u_1,),(\al+\be)(x_2 + u_2)].
\end{eqnarray}
The left hand is equal to
\begin{eqnarray*}
&&(\al+\be)([x_1, x_2] +\nu(x_1, x_2)+\rho(x_1)(u_2) -\rho(x_2)(u_1))\\
&=&\al([x_1, x_2]) +\be\circ\nu(x_1, x_2)+\be\circ\rho(x_1)(u_2) -\be\circ\rho(x_2)(u_1),
\end{eqnarray*}
and the right hand is equal to
\begin{eqnarray*}
&&[\al(x_1) + \be(u_1),\al(x_2) + \be(u_2)]\\
&=&[\al(x_1),\al(x_2)]+\nu(\al(x_1),\al(x_2)) +\rho(\al(x_1))\circ\be(u_2) -\rho(\al(x_2))\circ\be(u_1).
\end{eqnarray*}
Since $\al$ is an algebraic homomorphism and by conditions (HR01),(CC01), we obtain equality \eqref{semi01}.

The condition (HLY02) becomes
\begin{eqnarray}\label{semi02}
\notag&&(\al+\be)([x_1 + u_1, x_2 + u_2, x_3 + u_3])\\
&=&[(\al+\be)(x_1 + u_1), (\al+\be)(x_2 + u_2),(\al+\be)(x_3 + u_3)].
\end{eqnarray}
The left hand is equal to
\begin{eqnarray*}
&&(\al+\be)([x_1, x_2, x_3] + \DD(x_1, x_2)(u_3)-\theta(x_1, x_3)(u_2)+ \theta( x_2, x_3)(u_1))\\
&=&\al([x_1, x_2, x_3]) +\be\circ\omega(x_1, x_2, x_3)+\be\circ \DD(x_1, x_2)(u_3)-\be\circ\theta(x_1, x_3)(u_2)\\
&&+ \be\circ\theta( x_2, x_3)(u_1),
\end{eqnarray*}
and the right hand is equal to
\begin{eqnarray*}
&&[\al(x_1) + \be(u_1),\al(x_2) + \be(u_2),\al(x_3) + \be(u_3)]\\
&=&[\al(x_1),\al(x_2),\al(x_3)]  + \omega(\al(x_1), \al(x_2), \al(x_3))+\DD(\al(x_1),\al(x_2))\circ\be(u_3)\\
&&-\theta(\al(x_1),\al(x_3))\circ\be(u_2)+ \theta(\al(x_2),\al(x_3))\circ\be(u_1).
\end{eqnarray*}
Since $\al$ is an algebraic homomorphism  and by conditions (HR02), (CC02), we obtain equality \eqref{semi02}.


For (HLY1) and (HLY2), by definition we have
\begin{eqnarray*}
{[x_1 + u_1, x_1 + u_1]}&=&[x_1, x_1] +\rho(x_1)(u_1) -\rho(x_1)(u_1)=0,\\
\notag{[x_1 + u_1, x_1 + u_1, x_3 + u_3]}&=&[x_1, x_1, x_3] + \DD(x_1, x_1)(u_3)-\theta(x_1, x_3)(u_1)\\
&& + \theta(x_1, x_3)(u_1)=0.
\end{eqnarray*}


For (HLY3),  we have
\begin{eqnarray*}
&&[x_1 + u_1, x_2 + u_2, x_3 + u_3]_\omega+c.p.\\
&=&\{[x_1, x_2, x_3] +\underline{\omega(x_1,x_2,x_3)} +\DD(x_1, x_2)(u_3)-\theta(x_1, x_3)(u_2)\\
&& + \theta( x_2, x_3)(u_1)\}+c.p.
\end{eqnarray*}
and
\begin{eqnarray*}
&&[[x_1 + u_1, x_2 + u_2]_\nu, \al(x_3) + \be(u_3)]_\nu+c.p.\\
&=&[[x_1, x_2] +\rho(x_1)(u_2) -\rho(x_2)(u_1),\al(x_3) + \be(u_3)]_\nu+c.p.\\
&=&\{[[x_1, x_2], \al(x_3)] +\underline{\nu([x_1,x_2],x_3)} +\rho([x_1, x_2])\circ\be(u_3)\\
&& \underline{-\rho(\al(x_3))\nu(x_1, x_2)}-\rho(\al(x_3))\rho(x_1)(u_2) +\rho(\al(x_3))\rho(x_2)(u_1)\}+c.p.
\end{eqnarray*}
Thus by (HR31) and (CC1) we obtain
\begin{eqnarray}
[x_1 + u_1, x_2 + u_2, x_3 + u_3]_\omega+c.p.+[x_1 + u_1, x_2 + u_2]_\nu, \al(x_3) + \be(u_3)]_\nu+c.p.=0.
\end{eqnarray}

For (HLY4), we have
\begin{eqnarray*}
&&[[x_1 + u_1, x_2 + u_2]_\nu,\al(x_3) + \be(u_3),\al(y_1)+ \be(v_1)]_\omega+c.p.\\
&=&\{[[x_1, x_2], \al(x_3), \al(y_1)]  +\underline{\omega([x_1,x_2],\al(x_3),\al(y_1))} \\
&&+  \DD([x_1, x_2], \al(x_3))(\be(v_1))-\theta([x_1, x_2], \al(y_1))(\be(u_3)) \\
&&+ \underline{\theta(\al(x_3), \al(y_1))(\nu(x_1, x_2)}+\rho(x_1)(u_2) -\rho(x_2)(u_1))\}+c.p.\\
&=&0,
\end{eqnarray*}
where the last equality is by (CC2), (HR41) and (HR42).

For (HLY5), we have
\begin{eqnarray*}
&&[\al(x_1) + \be(u_1), \al(x_2) + \be(u_2), [y_1+v_1, y_2+v_2]_\nu]_\omega\\
&=&[\al(x_1), \al(x_2), [y_1, y_2]]+\underline{\omega(\al(x_1), \al(x_2), [y_1, y_2])}+\\
&&\underline{\DD(\al(x_1), \al(x_2))\Big(\nu(y_1, y_2)}+\rho(y_1)(v_2) -\rho(y_2)(v_1)\Big)\\
&&-\theta(x_1, [y_1, y_2])(u_2) + \theta(x_2, [y_1, y_2])(u_1),
\end{eqnarray*}
and
\begin{eqnarray*}
&&[[x_1 + u_1, x_2 + u_2, y_1+v_1]_\omega, \al^2(y_2)+\be^2(v_2)]_\nu\\
&&+[\al^2(y_1)+\be^2(v_1), [x_1 + u_1, x_2 + u_2, y_2+v_2]_\omega]_\nu\\
&=&[[x_1, x_2, y_1], \al^2(y_2)]+\underline{\nu([x_1, x_2, y_1], \al^2(y_2))}+ \rho([x_1, x_2, y_1])(v_2)\\
&&\underline{-\rho(y_2)\Big(\omega(x_1, x_2, y_1)}+\DD(x_1, x_2)(v_1)-\theta(x_1, y_1)(u_2) + \theta( x_2, y_1)(u_1)\Big)\\
&&+[\al^2(y_1), [x_1, x_2, y_2]]+\underline{\nu(\al^2(y_1), [x_1, x_2, y_2])} \\
&&+\underline{\rho(y_1)\Big(\omega(x_1, x_2, y_2)}+\DD(x_1, x_2)(v_2)-\theta(x_1, y_2)(u_2) + \theta( x_2, y_2)(u_1)\Big)\\
&& -\rho([x_1, x_2, y_2])(v_1).
\end{eqnarray*}
Thus by (CC3), (HR51) and (HR52) we obtain
\begin{eqnarray*}
&&[\al(x_1) + \be(u_1), \al(x_2) + \be(u_2), [y_1+v_1, y_2+v_2]_\nu]_\omega\\
&=&[[x_1 + u_1, x_2 + u_2, y_1+v_1]_\omega, \al^2(y_2)+\be^2(v_2)]_\nu\\
&&+[\al^2(y_1)+\be^2(v_1), [x_1 + u_1, x_2 + u_2, y_2+v_2]_\omega]_\nu
\end{eqnarray*}
 Therefore (HLY5) is valid.

Now it suffices to verify (HLY6). By definition,
\begin{eqnarray*}
&&[\al^2(x_1) + \be^2(u_1), \al^2(x_2) + \be^2(u_2),[y_1 + v_1, y_2 + v_2, y_3 + v_3]]\\
&=&[\al^2(x_1), \al^2(x_2), [y_1, y_2, y_3]] +\underline{\omega(\al^2(x_1), \al^2(x_2), [y_1, y_2, y_3])}\\
&&-\theta(\al^2(x_1), [y_1, y_2, y_3])(\be^2(u_2)) + \theta(\al^2(x_2), [y_1, y_2, y_3])(\be^2(u_1))\\
&& +\underline{\DD(\al^2(x_1), \al^2(x_2))\Big(\omega(y_1, y_2, y_3)}+\DD(y_1, y_2)(v_3)-\theta(y_1, y_3)(v_2) + \theta( y_2, y_3)(v_1)\Big),
\end{eqnarray*}
\begin{eqnarray*}
&&[[x_1 + u_1, x_2 + u_2,y_1 + v_1], \al^2(y_2) + \be^2(v_2), \al^2(y_3) + \be^2(v_3)]\\
&=&[[x_1, x_2, y_1], \al^2(y_2), \al^2(y_3)]]+\underline{\omega([x_1, x_2, y_1], \al^2(y_2), \al^2(y_3)])}\\
&& -\DD([x_1, x_2, y_1],\al^2(y_2))(\be^2(v_3)) + \theta([x_1, x_2, y_1],\al^2(y_3))(\be^2(u_1))\\
&& +\underline{\theta(\al^2(y_2), \al^2(y_3))\Big(\omega(x_1, x_2,y_1)}+\DD(x_1, x_2)(v_1)-\theta(x_1, y_1)(u_2) + \theta(x_2, y_1)(u_1)\Big),
\end{eqnarray*}
\begin{eqnarray*}
&&[\al^2(y_1) + \be^2(v_1), [x_1 + u_1, x_2 + u_2, y_2 + v_2],\al^2(y_3) + \be^2(v_3)]\\
&=&[\al^2(y_1), [x_1, x_2, y_2],\al^2(y_3)] +\omega(\al^2(y_1), [x_1, x_2, y_2],\al^2(y_3))\\
&&+\DD(\al^2(y_1), [x_1, x_2, y_2])(\be^2(v_3)) + \theta([x_1, x_2, y_2], \al^2(y_3))(\be^2(v_1))\\
&& \underline{-\theta(\al^2(y_1), \al^2(y_3))\Big(\omega(x_1, x_2,y_2)}+ \DD(x_1, x_2)(v_2)-\theta(x_1, y_2)(u_2) + \theta( x_2, y_2)(u_1)\Big),
\end{eqnarray*}
\begin{eqnarray*}
&&[\al^2(y_1) + \be^2(v_1), \al^2(y_2) + \be^2(v_2),[x_1 + u_1, x_2 + u_2, \al^2(y_3) + \be^2(v_3)]]\\
&=&[\al^2(y_1) + \al^2(y_2), [x_1, x_2, y_3]] +\underline{\omega(\al^2(y_1) + \al^2(y_2), [x_1, x_2, y_3])}\\
&&-\theta(\al^2(y_1), [x_1, x_2, y_3])(\be^2(v_2)) + \theta(\al^2(y_2), [x_1, x_2, y_3])(\be^2(v_1))\\
&& +\underline{\DD(\al^2(y_1), \al^2(y_2))\Big(\omega(x_1, x_2, y_3)}+\DD(x_1, x_2)(v_3)-\theta(x_1, y_3)(u_2) + \theta( x_2, y_3)(u_1)\Big),
\end{eqnarray*}
It follows that
\begin{eqnarray*}
&&[\al^2(x_1) + \be^2(u_1), \al^2(x_2) + \be^2(u_2),[y_1 + v_1, y_2 + v_2, y_3 + v_3]]\\
&=&[[x_1 + u_1, x_2 + u_2,y_1 + v_1], \al^2(y_2) + \be^2(v_2), \al^2(y_3) + \be^2(v_3)]\\
&&+[\al^2(y_1) + \be^2(v_1), [x_1 + u_1, x_2 + u_2, y_2 + v_2], y_3 + v_3]\\
&&+[\al^2(y_1) + \be^2(v_1), \al^2(y_2) + \be^2(v_2),[x_1 + u_1, x_2 + u_2, \al^2(y_3) + \be^2(v_3)]]
\end{eqnarray*}
by (CC4), (HR61) and (HR62).
Therefore we obtain a HLYA  on $T\oplus V$ under the maps \eqref{eq:newbracket00}, \eqref{eq:newbracket01} and \eqref{eq:newbracket02}.
The proof is completed.
\qed

\begin{lem}\label{thm:2-cocylce}
 Two abelian extensions of HLYAs $0\to\frkh{\to} T\oplus_{(\nu,\omega)}\h{\to} T\to 0$
 and  $0\to\frkh{\to} T\oplus_{(\nu',\omega')}\h{\to} T\to 0$ are equivalent if and only if $(\nu,\omega)$ and ${(\nu',\omega')}$ are in the same cohomology class.
\end{lem}

\pf Assume the two extensions are equivalent, we choose $F: T\oplus_{(\nu,\omega)}\h\to  T\oplus_{(\nu',\omega')}\h$ to be the corresponding homomorphism. Then we get
\begin{eqnarray}
\label{ggg}&&F[x_1, x_2]_{\nu}=[F(x_1),F(x_2)]_{\nu'},\\
\label{hhh}&&F[x_1, x_2, x_3]_{\omega}=[F(x_1),F(x_2),F(x_3)]_{\omega'}.
\end{eqnarray}
Since $F$ is an equivalence of extensions, there exist $f: T\to \frkh$ such that
$$F(x_i+u)=x_i+f(x_i)+u,\quad \forall x_i\in T.$$

Now the equation  \eqref{ggg} is equal to
\begin{eqnarray*}
&&[x_1, x_2]+f([x_1, x_2])+\nu(x_1, x_2)\\
&=&[x_1,x_2]+\nu'(x_1,x_2)+\rho(x_1)f(x_2)-\rho(x_2)f(x_1).
\end{eqnarray*}
Thus we have
\begin{eqnarray}\label{eq:exact4}
\nonumber(\nu-\nu')(x_1, x_2)&=&\rho(x_1)f(x_2)-\rho(x_2)f(x_1)-f([x_1, x_2]).
\end{eqnarray}

The equation \eqref{hhh} is equalent to
\begin{eqnarray*}
&&[x_1, x_2, x_3]+\omega(x_1, x_2, x_3)+f([x_1, x_2, x_3])\\
&=&[x_1,x_2,x_3]+\omega'(x_1,x_2,x_3)\\
&&+\DD(x_1, x_2)f(x_3)-\theta(x_1,x_3)f(x_2)+\theta(x_2,x_3)f(x_1).
\end{eqnarray*}
Thus we have
\begin{eqnarray}\label{eq:exact4}
\nonumber&&(\omega-\omega')(x_1, x_2, x_3)\\
&=&\DD(x_1, x_2)f(x_3)-\theta(x_1,x_3)f(x_2)+\theta(x_2,x_3)f(x_1)-f([x_1, x_2, x_3]).
\end{eqnarray}
Therefore $(\nu,\omega)$ and $(\nu',\omega')$ are in the same cohomology class.
Conversely, if $(\nu,\omega)$ and $(\nu',\omega')$  are in the same cohomology class, then we can show that $F$ is an equivalence.
We omit the details.
\qed
\medskip

Finally, we obtain the main result of this section:

\begin{thm}\label{thm:main99}
Let $T$ be a  HLYA and $V$ a $T$-module. Then there is a one-to-one correspondence between the set of equivalence classes of abelian extensions of the HLYA and the $(2,3)$-cohomology group. More precisely, there is a bijection map
$$\Ext( T,\frkh)\to H^2(T, V)\times H^3(T, V).$$
Therefore, the abelian extensions of $T$ by $V$ are  classified by the $(2,3)$-cohomology group.
\end{thm}

\emptycomment{
Let $T$ be a HLYA and $(V, \rho, D, \theta)$ a $T$-module.
For every linear map $f: T\to V$, according to Proposition \ref{prop:1cob}, the map $(\nu,\omega)$ defined by (BB1) and (BB2)
is a (2,3)-cocycle which is denoted by $(\nu,\omega)_f$.
Since this type of (2,3)-cocycles are in the same cohomology class with $(\nu',\omega')=(0,0)$, then by Theorem \ref{thm:main99}, we get

\begin{cor}
Let $T$ be a HLYA and $(V, \rho, D, \theta)$ a $T$-module.
Then for every linear map $f: T\to V$, the map $F: E_{(\nu,\omega)}\to E_{(\nu,\omega)+(\nu,\omega)_f}$ given by
$$F(x+u)=x+f(x)+u,\quad \forall x\in T, u\in V,$$
is a HLYA isomorphism.
\end{cor}

\begin{cor}
Let $T$ be a HLYA and $(V, \rho, D, \theta)$ a $T$-module.
\end{cor}
}

\section*{Acknowledgements}
The research was supported by Doctoral Research Program (5101019170129) of Henan Normal University.


\begin{thebibliography}{99}

\bibitem{Makhlouf01}
H. Ataguema, A. Makhlouf, and S. Silvestrov, Generalization of $n$-ary Nambu algebras and beyond, J. Math. Phys. 50(2009), 083501.

\bibitem{Makhlouf02}
J. Arnlind, A. Makhlouf, and S. Silvestrov, Ternary Hom-Nambu-Lie algebras induced by Hom-Lie algebras, J. Math. Phys. 51(2010), 043515.



\bibitem{Makhlouf03}
F. Ammar, S. Mabrouk, A. Makhlouf,
Representations and cohomology of $n$-ary multiplicative Hom-Nambu-Lie algebras, J. Geom. Physics,  61(2011), 1898--1913.



\bibitem{BEM01}
P. Benito, A. Elduque, and F. Mart\'{\i}n-Herce, Irreducible Lie-Yamaguti algebras, J. Pure Appl. Algebra 213(2009), 795--808.

\bibitem{Caenepeel}
S. Caenepeel and I. Goyvaerts, Monoidal Hom-Hopf Algebras,  Commu. Algebra 39(2011):  2216--2240.



\bibitem{Dor} I.~Dorfman, Dirac Structures and Integrability of Nonlinear Evolution Equation. John Wiley \& Sons, Ltd., Chichester, 1993.



\bibitem{Gerstenhaber1} M. Gerstenhaber, On the deformation of rings and algebras, Ann. Math. 79(1964): 59--103.

\bibitem{GNI}
D. Gaparayi and A. Nourou Issa,  A twisted Generalization of Lie-Yamaguti algebras,
Int. J. Algebra, 6(2012), 339--352.


\bibitem{Hom1}
J. T. Hartwig, D. Larsson, S. D. Silvestrov, Deformations of Lie algebras using $\sigma$-derivations, J. Algebra 295(2006), 314--361.


\bibitem{Kik01}
M. Kikkawa, Geometry of homogeneous Lie loops, Hiroshima Math. J. 5(1975), no. 2, 141--179.


\bibitem{KW}
M.K. Kinyon and A. Weinstein, Leibniz algebras, Courant algebroids, and multiplications on reductive homogeneous spaces,
 Amer. J. Math. 123(2001), no. 3, 525--550.


\bibitem{Ma}
Y. Ma, L.Y. Chen and J. Lin,
One-parameter formal deformations of Hom-Lie-Yamaguti algebras,
J. Math. Phys. 56(2015), 011701.



\bibitem{NR}
A. Nijenhuis and R. W. Richardson, Cohomology and deformations in graded Lie algebras,
Bull. Amer. Math. Soc. 72(1966), 1--29.


\bibitem{Nom}
K. Nomizu, Invariant affine connections on homogeneous spaces, Amer. J. Math 76(1954), 33--65.

\bibitem{Sheng}
Y. Sheng, Representations of hom-Lie algebras, Algebra and Representation Theory, 15(6)(2012), 1081--1098.

\bibitem{Ya1}
K. Yamaguti, On the Lie triple system and its generalization, J. Sci. of Hiroshima Univ., Ser. A, v. 21, 1958, pp. 155--160.


\bibitem{Ya2}
K. Yamaguti, On cohomology groups of general Lie triple systems, Kumamoto J. Sci., A 8(1969), 135--146.



\bibitem{Yau1}
D. Yau,  Hom-algebras and homology,  J. Lie Theory 19(2009), 409--421.

\bibitem{Yau2}
D. Yau, On $n$-ary Hom-Nambu and Hom-Nambu-Lie algebras, J. Geom. Phys. 62(2012), 506--522.

\bibitem{Zhang0}
T. Zhang, Notes on Cohomologies of Lie Triple Systems, J. Lie Theory, 24(4)(2014), 909--929.

\bibitem{Zhang1}
T. Zhang and J. Li, Deformations and extensions of Lie-Yamaguti algebras, to appear in Linear and Multilinear Algebra.
DOI:10.1080/03081087.2014.1000815


\end{thebibliography}
\end{document}